\documentclass[11pt]{amsart}
\usepackage{amsmath,amsfonts,amssymb}
\begin{document}

\title[Homogenization of the elliptic Dirichlet problem]{Homogenization of the elliptic Dirichlet problem:
operator error estimates in $L_2$}

\author{T.~A.~Suslina}

\thanks{Supported by RFBR (grant no. 11-01-00458-a)
and the Program of support of the leading scientific schools}


\address{St. Petersburg State University, Department of Physics, Ul'yanovskaya 3, Petrodvorets,
St.~Petersburg, 198504, Russia}

\email{suslina@list.ru}


\begin{abstract}
Let $\mathcal{O} \subset \mathbb{R}^d$ be a bounded domain of class $C^2$.
In the Hilbert space $L_2(\mathcal{O};\mathbb{C}^n)$,
we consider a matrix elliptic second order differential operator $\mathcal{A}_{D,\varepsilon}$ with
the Dirichlet boundary condition. Here $\varepsilon>0$ is the small parameter.
The coefficients of the operator are periodic and depend on
$\mathbf{x}/\varepsilon$. A sharp order operator error estimate
$\|\mathcal{A}_{D,\varepsilon}^{-1} - (\mathcal{A}_D^0)^{-1} \|_{L_2 \to L_2} \leq C \varepsilon$
is obtained. Here $\mathcal{A}^0_D$ is the effective operator with constant coefficients and
with the Dirichlet boundary condition.
\end{abstract}

\maketitle

\section*{Introduction}

The paper concerns homogenization theory of periodic differential operators (DO's).
A broad literature is devoted to homogenization problems in the small period limit.
First of all, we mention the books [BeLPa], [BaPan], [ZhKO].

\noindent\textbf{0.1. Operator-theoretic approach to homogenization problems.}
In a series of papers [BSu1-3] by M.~Sh.~Birman and T.~A.~Suslina
a new operator-theoretic (spectral) approach to homogenization problems was suggested and developed.
By this approach, the so-called operator error estimates in homogenization
problems for elliptic DO's were obtained.
Matrix elliptic DO's acting in $L_2(\mathbb{R}^d;\mathbb{C}^n)$ and admitting a factorization of the form
$\mathcal{A}_\varepsilon = b(\mathbf{D})^* g(\mathbf{x}/\varepsilon)b(\mathbf{D})$, $\varepsilon>0$,
were studied. Here $g(\mathbf{x})$ is a periodic matrix-valued function and $b(\mathbf{D})$
is a first order DO. The precise assumptions on $g(\mathbf{x})$ and $b(\mathbf{D})$
are described below in Section 1.

In [BSu1-3], the equation ${\mathcal{A}}_\varepsilon \mathbf{u}_\varepsilon + \mathbf{u}_\varepsilon
= \mathbf{F}$, where $\mathbf{F} \in L_2(\mathbb{R}^d;\mathbb{C}^n)$, was considered.
The behavior of the solution $\mathbf{u}_\varepsilon$ for small $\varepsilon$ was studied.
The solution $\mathbf{u}_\varepsilon$ converges in $L_2(\mathbb{R}^d;\mathbb{C}^n)$
to the solution $\mathbf{u}_0$ of the "homogenized"\ equation
${\mathcal{A}}^0 \mathbf{u}_0 + \mathbf{u}_0 = \mathbf{F}$, as $\varepsilon \to 0$.
Here ${\mathcal{A}}^0 = b(\mathbf{D})^* g^0 b(\mathbf{D})$ is the \textit{effective operator}
with the constant effective matrix $g^0$. In [BSu1], it was proved that
$$
\| \mathbf{u}_\varepsilon - \mathbf{u}_0 \|_{L_2(\mathbb{R}^d)}
\leq C \varepsilon \| \mathbf{F} \|_{L_2(\mathbb{R}^d)}.
$$
In operator terms it means that the resolvent
$({\mathcal{A}}_\varepsilon +I)^{-1}$ converges in the operator norm in $L_2(\mathbb{R}^d;\mathbb{C}^n)$
to the resolvent of the effective operator, as $\varepsilon \to 0$, and
$$
\| ({\mathcal{A}}_\varepsilon +I)^{-1} - ({\mathcal{A}}^0 +I)^{-1}
\|_{L_2(\mathbb{R}^d) \to L_2(\mathbb{R}^d)} \leq C \varepsilon.
\eqno(0.1)
$$

In [BSu2], more accurate approximation of the resolvent
$({\mathcal{A}}_\varepsilon +I)^{-1}$ in the operator norm in $L_2(\mathbb{R}^d;\mathbb{C}^n)$
with an error term $O(\varepsilon^2)$ was obtained.

In [BSu3], approximation of the resolvent $({\mathcal{A}}_\varepsilon +I)^{-1}$
in the norm of operators acting from $L_2(\mathbb{R}^d;\mathbb{C}^n)$ to the Sobolev space
$H^1(\mathbb{R}^d;\mathbb{C}^n)$ was found:
$$
\| ({\mathcal{A}}_\varepsilon +I)^{-1} - ({\mathcal{A}}^0 +I)^{-1} - \varepsilon K(\varepsilon) \|_{L_2(\mathbb{R}^d) \to H^1(\mathbb{R}^d)}
\le C \varepsilon;
\eqno(0.2)
$$
this corresponds to approximation of $\mathbf{u}_\varepsilon$ in the "energy"\ norm.
Here $K(\varepsilon)$ is a corrector. It contains rapidly oscillating factors and so depends on
$\varepsilon$.

Estimates (0.1), (0.2) are called the \textit{operator error estimates}. They are order-sharp;
the constants in estimates are controlled explicitly in terms of the problem data.
The method of [BSu1--3] is based on the scaling transformation, the Floquet-Bloch theory and
the analytic perturbation theory.

\smallskip\noindent\textbf{0.2. A different approach} to operator error estimates in homogenization
problems was suggested by V.~V.~Zhikov. In [Zh1, Zh2, ZhPas, Pas], the scalar elliptic operator
$- \text{div}\, g(\mathbf{x}/\varepsilon) \nabla$ (where $g(\mathbf{x})$ is a matrix with
real entries) and the system of elasticity theory were studied.
Estimates of the form (0.1), (0.2) for the corresponding problems in $\mathbb{R}^d$ were obtained.
The method was based on analysis of the first order approximation to the solution and introducing
of an additional parameter. Besides the problems in $\mathbb{R}^d$, homogenization problems in
a bounded domain $\mathcal{O} \subset \mathbb{R}^d$ with the Dirichlet or Neumann boundary condition
were studied. Approximation of the solution in $H^1(\mathcal{O})$ was deduced from
the corresponding result in $\mathbb{R}^d$. Due to the "boundary layer"\ influence,
estimates in a bounded domain become worse and the error term is $O(\varepsilon^{1/2})$.
The estimate $\|\mathbf{u}_\varepsilon - \mathbf{u}_0\|_{L_2(\mathcal{O})}
\leq C \varepsilon^{1/2}\|\mathbf{F}\|_{L_2(\mathcal{O})}$
follows from approximation of the solution in $H^1(\mathcal{O})$ by roughening.

Similar results for the operator $- \text{div}\, g(\mathbf{x}/\varepsilon) \nabla$ in a bounded domain
with the Dirichlet or Neumann boundary condition were obtained
in the papers [Gr1, Gr2] by G.~Griso by the "unfolding"\ method.

\smallskip\noindent\textbf{0.3. Approximation of the resolvent in the $(L_2 \to H^1)$-norm.}
The present paper relies on the results of [PSu].
In that paper, matrix DO's $\mathcal{A}_{D,\varepsilon}$ in a bounded domain
$\mathcal{O}\subset \mathbb{R}^d$ of class $C^2$ were studied.
The operator $\mathcal{A}_{D,\varepsilon}$ is defined by
the differential expression $b(\mathbf{D})^* g(\mathbf{x}/\varepsilon) b(\mathbf{D})$
with the Dirichlet condition on $\partial \mathcal{O}$.
The effective operator $\mathcal{A}_D^0$ is given by the expression
$b(\mathbf{D})^* g^0 b(\mathbf{D})$ with the Dirichlet boundary condition.
The behavior for small $\varepsilon$ of the solution $\mathbf{u}_\varepsilon$ of the equation
$\mathcal{A}_{D,\varepsilon} \mathbf{u}_\varepsilon =\mathbf{F}$, where
$\mathbf{F} \in L_2(\mathcal{O};\mathbb{C}^n)$, is studied.
Estimates for the $H^1$-norm of the difference of the solution $\mathbf{u}_\varepsilon$
and its first order approximation are obtained.
By roughening of this result, an estimate for
$\|\mathbf{u}_\varepsilon - \mathbf{u}_0\|_{L_2(\mathcal{O})}$ is proved.
Here $\mathbf{u}_0$ is the solution of the equation $\mathcal{A}^0_D \mathbf{u}_0 =\mathbf{F}$.

In operator terms, the following estimates are obtained:
$$
\| \mathcal{A}_{D,\varepsilon}^{-1} - (\mathcal{A}^0_D)^{-1} - \varepsilon K_D(\varepsilon) \|_{L_2(\mathcal{O}) \to H^1(\mathcal{O})}
\leq C \varepsilon^{1/2},
\eqno(0.3)
$$
$$
\| \mathcal{A}_{D,\varepsilon}^{-1} - (\mathcal{A}^0_D)^{-1} \|_{L_2(\mathcal{O}) \to L_2(\mathcal{O})}
\leq C \varepsilon^{1/2}.
\eqno(0.4)
$$
Here $K_D(\varepsilon)$ is the corresponding corrector.

The method of [PSu] is based on using estimates (0.1), (0.2)
for homogenization problem in $\mathbb{R}^d$ obtained in [BSu1,3]
and on the tricks suggested in [Zh2], [ZhPas] that allow one to deduce estimate (0.3) from (0.1), (0.2).
Main difficulties are related to estimating of the "discrepancy"\ $\mathbf{w}_\varepsilon$, which satisfies
the equation $\mathcal{A}_\varepsilon \mathbf{w}_\varepsilon =0$ in $\mathcal{O}$ and the boundary
condition $\mathbf{w}_\varepsilon = \varepsilon K_D(\varepsilon)\mathbf{F}$ on $\partial \mathcal{O}$.

\smallskip\noindent\textbf{0.4. The main result.}
It must be mentioned that estimate (0.4) is quite a rough consequence of (0.3).
So, the refinement of estimate (0.4) is a natural problem.
In [ZhPas], for the case of the scalar elliptic operator
$- \text{div}\, g(\mathbf{x}/\varepsilon) \nabla$ (where $g(\mathbf{x})$ is a matrix with real entries)
an estimate for
$\|\mathcal{A}_{D,\varepsilon}^{-1} - (\mathcal{A}^0_D)^{-1} \|_{L_2\to L_2}$
of order $\varepsilon^{\frac{d}{2d-2}}$ for $d \geq 3$ and of order $\varepsilon |\log \varepsilon|$
for $d=2$ was obtained. The proof essentially relies on using the maximum principle which is
specific for scalar elliptic equations.

In the present paper, we prove a \textit{sharp order operator error estimate}
$$
\| \mathcal{A}_{D,\varepsilon}^{-1} - (\mathcal{A}^0_D)^{-1} \|_{L_2(\mathcal{O}) \to L_2(\mathcal{O})}
\leq C \varepsilon.
\eqno(0.5)
$$
Estimate (0.5) for matrix elliptic DO's refines even
the known classical (non-operator) error estimates.

Method of the proof relies on the results and technique of [PSu].
The problem reduces to estimating of the $L_2$-norm of $\mathbf{w}_\varepsilon$.
Using of the operator approach and duality arguments is important.
Employing approximation of the resolvent $\mathcal{A}_{D,\varepsilon}^{-1}$
in the norm of operators acting from $L_2(\mathcal{O};\mathbb{C}^n)$ to $H^1_0(\mathcal{O};\mathbb{C}^n)$,
we find approximation of the same operator in
the norm of operators acting from $H^{-1}(\mathcal{O};\mathbb{C}^n)$ to $L_2(\mathcal{O};\mathbb{C}^n)$.
The last approximation combined with the boundary layer estimates allows one to obtain
the required estimate for the $L_2$-norm of $\mathbf{w}_\varepsilon$.

\smallskip\noindent\textbf{0.5. The plan of the paper.} The paper contains three sections.
In Section 1, the class of operators is introduced, the effective operator is described, and
the main result is formulated.
Section 2 contains some auxiliary statements needed for further investigation.
In Section 3, the main result is proved.

\smallskip\noindent\textbf{0.6. Notation.} Let $\mathfrak{H}$ and $\mathfrak{H}_*$ be complex
separable Hilbert spaces. The symbols $(\cdot,\cdot)_{\mathfrak{H}}$ and $\|\cdot\|_{\mathfrak{H}}$
stand for the inner product and the norm in $\mathfrak{H}$;
the symbol $\|\cdot\|_{\mathfrak{H} \to \mathfrak{H}_*}$
denotes the norm of a linear continuous operator acting from $\mathfrak{H}$ to $\mathfrak{H}_*$.

The symbols $\langle \cdot, \cdot \rangle$ and $|\cdot|$ stand for the inner product
and the norm in $\mathbb{C}^n$; $\mathbf{1} = \mathbf{1}_n$ is the identity $(n\times n)$-matrix.
We use the notation $\mathbf{x} = (x_1,\dots,x_d)\in \mathbb{R}^d$, $iD_j = \partial_j = \partial/\partial x_j$,
$j=1,\dots,d$, $\mathbf{D} = -i \nabla = (D_1,\dots,D_d)$.
The $L_p$-classes of $\mathbb{C}^n$-valued functions in a domain ${\mathcal O} \subset \mathbb{R}^d$
are denoted by $L_p({\mathcal O};\mathbb{C}^n)$, $1 \le p \leq \infty$.
The Sobolev classes of $\mathbb{C}^n$-valued functions in a domain ${\mathcal O} \subset \mathbb{R}^d$
are denoted by $H^s({\mathcal O};\mathbb{C}^n)$.
By $H^1_0(\mathcal{O};\mathbb{C}^n)$ we denote the closure of $C_0^\infty(\mathcal{O};\mathbb{C}^n)$
in $H^1(\mathcal{O};\mathbb{C}^n)$.
If $n=1$, we write simply $L_p({\mathcal O})$, $H^s({\mathcal O})$, etc., but sometimes
we use such abbreviated notation also for spaces of vector-valued or matrix-valued functions.

\smallskip\noindent\textbf{0.7. Acknowledgement.} The author is grateful to A.~I.~Nazarov for
fruitful stimulating discussions.

\section*{\S 1. Statement of the problem. Results}

\smallskip\noindent\textbf{1.1. The class of operators.}
Let $\Gamma \subset \mathbb{R}^d$ be a lattice, and let $\Omega \subset \mathbb{R}^d$
be the elementary cell of the lattice $\Gamma$.
We denote $|\Omega| = \text{meas}\, \Omega$.
Below $\widetilde{H}^1(\Omega)$ stands for the subspace of
functions in $H^1(\Omega)$ whose $\Gamma$-periodic extension to
$\mathbb{R}^d$ belongs to $H^1_{\text{loc}}(\mathbb{R}^d)$.
If $\varphi(\mathbf{x})$ is a $\Gamma$-periodic function in $\mathbb{R}^d$, we denote
$$
\varphi^\varepsilon(\mathbf{x}) := \varphi(\varepsilon^{-1}\mathbf{x}),\quad \varepsilon >0.
$$

Let $\mathcal{O} \subset \mathbb{R}^d$ be a bounded domain of class $C^2$.
In $L_2(\mathcal{O};\mathbb{C}^n)$, we define an operator $\mathcal{A}_{D,\varepsilon}$
formally given by the differential expression
$$
\mathcal{A}_\varepsilon = b(\mathbf{D})^* g^\varepsilon(\mathbf{x}) b (\mathbf{D})
\eqno(1.1)
$$
with the Dirichlet condition on $\partial \mathcal{O}$.
Here $g(\mathbf{x})$ is a measurable $(m \times m)$-matrix-valued
function (in general, with complex entries). We assume that $g(\mathbf{x})$ is periodic with respect
to the lattice $\Gamma$, bounded and uniformly positive definite.
Next, $b(\mathbf{D}) = \sum_{l=1}^d b_l D_l$ is an $(m\times n)$-matrix first order DO
with constant coefficients. Here $b_l$ are constant matrices (in general, with complex entries).
The symbol $b(\boldsymbol{\xi}) = \sum_{l=1}^d b_l \xi_l$, $\boldsymbol{\xi} \in \mathbb{R}^d$,
corresponds to the operator $b(\mathbf{D})$.
It is assumed that $m \ge n$ and that
$\textrm{rank} \, b (\boldsymbol{\xi}) = n, \ \forall \boldsymbol{\xi} \neq 0$.
This condition is equivalent to the following inequalities
$$
\alpha_0 \mathbf{1}_n \leq b(\boldsymbol{\theta})^* b(\boldsymbol{\theta}) \leq \alpha_1 \mathbf{1}_n,
\quad \boldsymbol{\theta} \in \mathbb{S}^{d-1}, \quad 0 < \alpha_0 \leq \alpha_1 < \infty,
\eqno(1.2)
$$
with some positive constants $\alpha_0$ and $\alpha_1$.

The precise definition is the following: $\mathcal{A}_{D,\varepsilon}$ is the selfadjoint operator in
$L_2(\mathcal{O};\mathbb{C}^n)$ generated by the quadratic form
$$
a_{D,\varepsilon}[\mathbf{u},\mathbf{u}] = \int_{\mathcal{O}} \left\langle g^\varepsilon (\mathbf{x})
b (\mathbf{D}) \mathbf{u}, b(\mathbf{D}) \mathbf{u} \right\rangle \,
d \mathbf{x}, \quad \mathbf{u} \in H^1_0 (\mathcal{O}; \mathbb{C}^{n}).
$$
Under the above assumptions this form is closed in $L_2(\mathcal{O};\mathbb{C}^n)$
and positive definite. Moreover, we have
$$
c_0 \int_{\mathcal{O}} |\mathbf{D} \mathbf{u}|^2\, d\mathbf{x} \le a_{D,\varepsilon}[\mathbf{u}, \mathbf{u}] \le c_1 \int_{\mathcal{O}} |\mathbf{D} \mathbf{u}|^2\, d\mathbf{x},
\ \ \mathbf{u} \in H^1_0(\mathcal{O};\mathbb{C}^n),
\eqno(1.3)
$$
where $c_0 = \alpha_0 \|g^{-1}\|^{-1}_{L_\infty}$, $c_1 = \alpha_1 \|g\|_{L_\infty}$.
It is easy to check (1.3) extending
$\mathbf{u}$ by zero to $\mathbb{R}^d\setminus \mathcal{O}$, using
the Fourier transformation and taking (1.2) into account.

The simplest example of the operator (1.1) is the scalar elliptic operator
$\mathcal{A}_\varepsilon = -\text{div}\, g^\varepsilon(\mathbf{x}) \nabla =
\mathbf{D}^* g^\varepsilon(\mathbf{x})\mathbf{D}.$
In this case we have $n=1$, $m=d$, $b(\mathbf{D})=\mathbf{D}$. Obviously, condition (1.2)
is valid with $\alpha_0=\alpha_1=1$.
Another example is the operator of elasticity theory which can be written in the form (1.1)
with $n=d$, $m=d(d+1)/2$. These and other examples are discussed in [BSu1] in detail.

\textit{Our goal} is to find approximation for small $\varepsilon$
for the operator $\mathcal{A}_{D,\varepsilon}^{-1}$ in the operator norm in $L_2(\mathcal{O};\mathbb{C}^n)$.
In terms of solutions, we are interested in the behavior of the generalized solution
$\mathbf{u}_\varepsilon \in H^1_0(\mathcal{O};\mathbb{C}^n)$ of the Dirichlet problem
$$
b (\mathbf{D})^* g^{\varepsilon} (\mathbf{x}) b (\mathbf{D})
\mathbf{u}_{\varepsilon}(\mathbf{x}) = \mathbf{F}(\mathbf{x}),\ \ \mathbf{x} \in \mathcal{O};
\quad \mathbf{u}_\varepsilon\vert_{\partial \mathcal{O}}=0,
\eqno(1.4)
$$
where $\mathbf{F} \in L_2 (\mathcal{O}; \mathbb{C}^{n})$.
Then $\mathbf{u}_\varepsilon = \mathcal{A}_{D,\varepsilon}^{-1} \mathbf{F}$.

\smallskip\noindent\textbf{1.2. The effective operator.}
In order to formulate the results, we need to introduce the effective operator $\mathcal{A}^0_D$.

Let an $(n\times m)$-matrix-valued function $\Lambda(\mathbf{x})$ be the (weak)
$\Gamma$-periodic solution of the problem
$$
b(\mathbf{D})^* g (\mathbf{x})\left( b(\mathbf{D}) \Lambda(\mathbf{x}) + \mathbf{1}_m \right) = 0,
\quad \int_{\Omega} \Lambda(\mathbf{x}) \, d \mathbf{x} = 0.
\eqno(1.5)
$$
In other words, for the columns $\mathbf{v}_j(\mathbf{x})$, $j=1,\dots,m,$
of the matrix $\Lambda(\mathbf{x})$ the following is true:
$\mathbf{v}_j \in \widetilde{H}^1(\Omega;\mathbb{C}^n)$, we have
$$
\int_\Omega \langle g(\mathbf{x}) (b(\mathbf{D})\mathbf{v}_j(\mathbf{x}) + \mathbf{e}_j),
b(\mathbf{D}) \boldsymbol{\eta}(\mathbf{x})\rangle\,d\mathbf{x} =0,
\quad \forall \boldsymbol{\eta} \in \widetilde{H}^1(\Omega;\mathbb{C}^n),
$$
and $\int_\Omega \mathbf{v}_j(\mathbf{x})\,d\mathbf{x}=0$.
Here $\mathbf{e}_1,\dots,\mathbf{e}_m$ is the standard orthonormal basis in $\mathbb{C}^m$.

The so-called \textit{effective matrix} $g^0$ of size $m\times m$ is defined as follows:
$$
g^0 = |\Omega|^{-1}
\int_{\Omega} g (\mathbf{x}) \left( b (\mathbf{D})
\Lambda (\mathbf{x}) + \mathbf{1}_m \right)\, d \mathbf{x}.
\eqno(1.6)
$$
It turns out that the matrix (1.6) is positive definite.
The \textit{effective operator} $\mathcal{A}^0_D$ for $\mathcal{A}_{D,\varepsilon}$
is given by the differential expression
$$
\mathcal{A}^0 = b (\mathbf{D})^* g^0 b (\mathbf{D})
$$
with the Dirichlet condition on $\partial \mathcal{O}$.
The domain of this operator is $H^1_0(\mathcal{O};\mathbb{C}^n) \cap H^2(\mathcal{O};\mathbb{C}^n)$
(see Subsection 2.2 below).

Consider the "homogenized"\ Dirichlet problem
$$
b(\mathbf{D})^* g^0 b(\mathbf{D}) \mathbf{u}_0(\mathbf{x}) = \mathbf{F}(\mathbf{x}),
\ \ \mathbf{x} \in \mathcal{O};
\quad \mathbf{u}_0\vert_{\partial \mathcal{O}}=0.
\eqno(1.7)
$$
Then $\mathbf{u}_0 = (\mathcal{A}^0_D)^{-1}\mathbf{F}$.
As $\varepsilon \to 0$, the solution $\mathbf{u}_\varepsilon$ of the problem (1.4) converges in
$L_2(\mathcal{O};\mathbb{C}^n)$ to $\mathbf{u}_0$; for operators of the form (1.1) this was proved in [PSu].
We wish to estimate $\|\mathbf{u}_\varepsilon - \mathbf{u}_0\|_{L_2(\mathcal{O})}$.

\smallskip\noindent\textbf{1.3. The main result}. Denote
$$
(\partial \mathcal{O})_\varepsilon =
\{\mathbf{x} \in \mathbb{R}^d:\ \text{dist}\,\{\mathbf{x},\partial \mathcal{O}\} < \varepsilon\}.
$$

Now we formulate the main result.

\smallskip\noindent\textbf{Theorem 1.1.} \textit{Assume that} $\mathcal{O}\subset \mathbb{R}^d$
\textit{is a bounded domain of class} $C^2$. \textit{Let} $g(\mathbf{x})$
\textit{and} $b(\mathbf{D})$ \textit{satisfy the assumptions of Subsection} 1.1.
\textit{Let} $\mathbf{u}_{\varepsilon}$ \textit{be the solution of the problem} (1.4), \textit{and let}
$\mathbf{u}_0$ \textit{be the solution of the problem} (1.7) \textit{with}
$\mathbf{F} \in L_2(\mathcal{O};\mathbb{C}^n)$.
\textit{Let} $\varepsilon_1 \in (0,1]$
\textit{be such that the set} $(\partial \mathcal{O})_{\varepsilon_1}$
\textit{can be covered by a finite number of open sets admitting diffeomorphisms of class} $C^2$
\textit{rectifying the boundary} $\partial \mathcal{O}$. \textit{Let} $2r_1 = \text{diam}\,\Omega$,
$\varepsilon_2 = \varepsilon_1(1+r_1)^{-1}$, \textit{and} $\varepsilon_0 = \varepsilon_2/2$.
\textit{Then for} $0< \varepsilon \leq \varepsilon_0$ \textit{we have}
$$
\| \mathbf{u}_{\varepsilon} - \mathbf{u}_0
\|_{L_2(\mathcal{O}; \mathbb{C}^{n})} \leq C_1 \varepsilon
\| \mathbf{F} \|_{L_2 (\mathcal{O}; \mathbb{C}^{n})},
\eqno(1.8)
$$
\textit{or, in operator terms,}
$$
\| \mathcal{A}_{D,\varepsilon}^{-1} - (\mathcal{A}^0_D)^{-1}
\|_{L_2 (\mathcal{O}; \mathbb{C}^{n}) \to
L_2 (\mathcal{O}; \mathbb{C}^{n})} \leq C_1 \varepsilon.
$$
\textit{The constant} $C_1$ \textit{depends only on}
$m$, $d$, $\alpha_0$, $\alpha_1$, $\| g \|_{L_\infty}$,
$\| g^{-1} \|_{L_\infty}$, \textit{the parameters of the lattice} $\Gamma$,
\textit{and the domain} $\mathcal{O}$.

\section*{\S 2. Auxiliary statements}

\smallskip\noindent\textbf{2.1. The energy inequality.}
Consider the problem (1.4) with the right-hand side of class $H^{-1}(\mathcal{O};\mathbb{C}^n)$.
Recall that $H^{-1}(\mathcal{O};\mathbb{C}^n)$ is defined as the space dual to
$H^1_0(\mathcal{O};\mathbb{C}^n)$ with respect to the
$L_2(\mathcal{O};\mathbb{C}^n)$-coupling. If $\mathbf{f} \in H^{-1}(\mathcal{O};\mathbb{C}^n)$ and
$\boldsymbol{\eta} \in H^1_0 (\mathcal{O}; \mathbb{C}^{n})$, then the symbol
$(\mathbf{f}, \boldsymbol{\eta})_{L_2(\mathcal{O})}=\int_\mathcal{O}
\langle \mathbf{f}, \boldsymbol{\eta} \rangle\,d\mathbf{x}$
stands for the value of the functional $\mathbf{f}$ on the element $\boldsymbol{\eta}$. Herewith,
$$
\left| \int_\mathcal{O} \langle \mathbf{f}, \boldsymbol{\eta} \rangle\,d\mathbf{x}\right|
\leq \|\mathbf{f}\|_{H^{-1}(\mathcal{O};\mathbb{C}^n)} \|\boldsymbol{\eta}\|_{H^1(\mathcal{O};\mathbb{C}^n)}.
$$
The following (standard) statement was checked in [PSu, Lemma 4.1].

\smallskip\noindent\textbf{Lemma 2.1.} \textit{Let} $\mathbf{f} \in H^{-1}(\mathcal{O};\mathbb{C}^n)$.
\textit{Suppose that} $\mathbf{z}_\varepsilon \in H^1_0(\mathcal{O};\mathbb{C}^n)$
\textit{is the generalized solution of the Dirichlet problem}
$$
b (\mathbf{D})^* g^{\varepsilon} (\mathbf{x}) b (\mathbf{D})
\mathbf{z}_{\varepsilon}(\mathbf{x}) = \mathbf{f}(\mathbf{x}),\ \ \mathbf{x} \in \mathcal{O};
\quad \mathbf{z}_\varepsilon\vert_{\partial \mathcal{O}}=0.
$$
\textit{In other words}, $\mathbf{z}_\varepsilon$ \textit{satisfies the identity}
$$
\int_{\mathcal{O}} \langle g^{\varepsilon}(\mathbf{x}) b (\mathbf{D}) \mathbf{z}_{\varepsilon},
b (\mathbf{D}) \boldsymbol{\eta} \rangle\,d\mathbf{x} = \int_\mathcal{O} \langle \mathbf{f},
\boldsymbol{\eta} \rangle
\,d\mathbf{x}, \quad \forall \ \boldsymbol{\eta} \in H^1_0 (\mathcal{O}; \mathbb{C}^{n}).
$$
\textit{Then the following "energy inequality"\  is true}:
$$
\| \mathbf{z}_{\varepsilon} \|_{H^1 (\mathcal{O}; \mathbb{C}^{n})} \leq
\widehat{C} \| \mathbf{f} \|_{H^{-1} (\mathcal{O}; \mathbb{C}^{n})},
$$
\textit{where} $\widehat{C}= (1 + (\text{diam}\,{\mathcal{O}})^2) \alpha_0^{-1} \| g^{-1} \|_{L_\infty}$.

\smallskip
It follows from Lemma 2.1 that the operator $\mathcal{A}_{D,\varepsilon}^{-1}$
acting in $L_2(\mathcal{O};\mathbb{C}^n)$ can be extended to a linear continuous operator
acting from $H^{-1}(\mathcal{O};\mathbb{C}^n)$ to $H^1_0(\mathcal{O};\mathbb{C}^n)$.
Applying Lemma 2.1 with $g^\varepsilon$ replaced by $g^0$, we see that the same
statement is true for the operator $(\mathcal{A}^0_D)^{-1}$.

Note that
$$
( \mathcal{A}^{-1}_{D,\varepsilon} \mathbf{f}_1, \mathbf{f}_2)_{L_2(\mathcal{O})}=
(\mathbf{f}_1, \mathcal{A}^{-1}_{D,\varepsilon}\mathbf{f}_2)_{L_2(\mathcal{O})},\quad
\mathbf{f}_1, \mathbf{f}_2 \in H^{-1}(\mathcal{O};\mathbb{C}^n).
\eqno(2.1)
$$
A similar identity is valid for the operator $(\mathcal{A}_D^0)^{-1}$.

All the statements of Subsection 2.1 are valid in arbitrary bounded domain $\mathcal{O}$
(without assumption that $\partial \mathcal{O} \in C^2$).

\smallskip\noindent\textbf{2.2. Properties of the solution of the homogenized problem.}
Due to the assumption $\partial \mathcal{O} \in C^2$, the solution $\mathbf{u}_0$ of the problem
(1.7) satisfies
$
\mathbf{u}_0 \in H^1_0 (\mathcal{O}; \mathbb{C}^{n}) \cap H^2 (\mathcal{O}; \mathbb{C}^{n}),
$
and
$$
\| \mathbf{u}_0 \|_{H^2 (\mathcal{O}; \mathbb{C}^{n})} \leq \widehat{c}
\| \mathbf{F} \|_{L_2 (\mathcal{O}; \mathbb{C}^{n})}.
\eqno(2.2)
$$
In operator terms, it means that
$$
\| (\mathcal{A}^0_D)^{-1}\|_{L_2(\mathcal{O};\mathbb{C}^n) \to H^2(\mathcal{O};\mathbb{C}^n)} \le \widehat{c}.
\eqno(2.3)
$$
The constant $\widehat{c}$ depends only on $\alpha_0$, $\alpha_1$,
$\|g\|_{L_\infty}$, $\|g^{-1}\|_{L_\infty}$, and the domain $\mathcal{O}$.
To justify these properties, it suffices to note that the operator $b(\mathbf{D})^* g^0 b(\mathbf{D})$ is
a \textit{strongly elliptic} matrix DO
and to apply the "additional smoothness"\ theorems for solutions of strongly elliptic systems
(see, e.~g., [McL, Chapter 4]).

\smallskip\noindent\textbf{2.3. Trace lemma.} We need the following simple statement;
see, e.~g., [PSu, Lemma 5.1].

\smallskip\noindent\textbf{Lemma 2.2.}
\textit{Denote}
$B_{\varepsilon} = \left\{ \mathbf{x} \in \mathcal{O}: \textrm{\upshape dist\itshape}
\,\{ \mathbf{x}, \partial \mathcal{O} \} < \varepsilon \right\}$.
\textit{Then for any} $z \in H^1 (\mathcal{O})$ \textit{we have}
$$
\int_{B_{\varepsilon}} |z|^{2} d \mathbf{x}
\leq \beta \varepsilon \| z \|_{H^1 (\mathcal{O})}\| z \|_{L_2 (\mathcal{O})},
\quad 0< \varepsilon \leq \varepsilon_1.
$$
\textit{Here} $\varepsilon_1$ \textit{is the same as in Theorem} 1.1.
\textit{The constant} $\beta$ \textit{depends only on the domain} $\mathcal{O}$.

\smallskip
Note that the statement of Lemma 2.2 is valid for any bounded domain $\mathcal{O}$ of class $C^1$.

\smallskip\noindent\textbf{2.4. Smoothing in Steklov's sense.}
Let $S_\varepsilon$ be the operator in $L_2(\mathbb{R}^d;\mathbb{C}^m)$ given by
$$
(S_\varepsilon \mathbf{u})(\mathbf{x}) = |\Omega|^{-1} \int_\Omega
\mathbf{u}(\mathbf{x} - \varepsilon \mathbf{z})\, d\mathbf{z}.
\eqno(2.4)
$$
It is said that the operator $S_\varepsilon$ is \textit{smoothing in Steklov's sense}.

We need the following property of the operator (2.4) (see [ZhPas, Lemma 1.1] or [PSu, Proposition 3.2]).

\smallskip\noindent\textbf{Lemma 2.3.} \textit{Let} $f(\mathbf{x})$
\textit{be a} $\Gamma$-\textit{periodic function in} $\mathbb{R}^d$ \textit{such that}
$f \in L_2(\Omega)$. \textit{Let} $[f^\varepsilon]$ \textit{denote the operator of multiplication
by the function} $f^\varepsilon(\mathbf{x})$. \textit{Then the operator} $[f^\varepsilon]S_\varepsilon$
\textit{is continuous in} $L_2(\mathbb{R}^d;\mathbb{C}^m)$, \textit{and}
$$
\| [f^\varepsilon]S_\varepsilon \|_{L_2(\mathbb{R}^d;\mathbb{C}^m)\to L_2(\mathbb{R}^d;\mathbb{C}^m)}
\leq |\Omega|^{-1/2} \| f \|_{L_2(\Omega)}.
$$

\smallskip\noindent\textbf{2.5. Properties of the matrix} $\Lambda(\mathbf{x})$.
Let $\widetilde{\Gamma}$ be the lattice dual to $\Gamma$.
By $\widetilde{\Omega}$ we denote the central Brillouin zone of $\widetilde{\Gamma}$, i.~e.,
$\widetilde{\Omega} = \{ \mathbf{k} \in \mathbb{R}^d:\
|\mathbf{k}| < |\mathbf{k}- \mathbf{b}|,\ 0\ne \mathbf{b} \in \widetilde{\Gamma}\}$.
Let $r_0$ be the radius of the ball inscribed in $\text{clos}\, \widetilde{\Omega}$.

Recall that the matrix-valued function $\Lambda(\mathbf{x})$ is the
$\Gamma$-periodic solution of the problem (1.5). In [BSu2, Subsection 7.3] it was proved that
$$
 \|\Lambda\|_{L_2(\Omega)}
\leq m^{1/2} (2r_0)^{-1} |\Omega|^{1/2}\alpha_0^{-1/2}
\|g\|^{1/2}_{L_\infty} \|g^{-1}\|^{1/2}_{L_\infty}.
\eqno(2.5)
$$

\smallskip
Let $[\Lambda^\varepsilon]$ be the operator of multiplication by the matrix-valued function
$\Lambda^\varepsilon(\mathbf{x})$; this operator acts from $L_2(\mathbb{R}^d;\mathbb{C}^m)$ to
$L_2(\mathbb{R}^d;\mathbb{C}^n)$. By Lemma 2.3 and estimate (2.5),
the norm of the operator $[\Lambda^\varepsilon] S_\varepsilon$ satisfies the following estimate:
$$
\begin{aligned}
&\|[\Lambda^\varepsilon] S_\varepsilon\|_{L_2(\mathbb{R}^d;\mathbb{C}^m)\to L_2(\mathbb{R}^d;\mathbb{C}^n)}
\leq |\Omega|^{-1/2} \|\Lambda\|_{L_2(\Omega)}
\\
&\leq m^{1/2} (2r_0)^{-1} \alpha_0^{-1/2} \|g\|^{1/2}_{L_\infty} \|g^{-1}\|^{1/2}_{L_\infty}=:M.
\end{aligned}
\eqno(2.6)
$$

\section*{\S 3. Proof of Theorem 1.1}

The proof of Theorem 1.1 relies on the results of [PSu],
where approximation of $\mathcal{A}_{D,\varepsilon}^{-1}$ in the norm of operators acting
from $L_2(\mathcal{O};\mathbb{C}^n)$ to $H^1(\mathcal{O};\mathbb{C}^n)$ was obtained.

\smallskip\noindent\textbf{3.1. Error estimates in} $H^1$.
We fix a linear continuous extension operator
$$
P_\mathcal{O}: H^2(\mathcal{O};\mathbb{C}^n) \to H^2(\mathbb{R}^d;\mathbb{C}^n)
\eqno(3.1)
$$
and put $\widetilde{\mathbf{u}}_0 = P_\mathcal{O} \mathbf{u}_0$. Then
$$
\|\widetilde{\mathbf{u}}_0 \|_{H^2(\mathbb{R}^d;\mathbb{C}^n)}
\leq C_\mathcal{O} \| \mathbf{u}_0 \|_{H^2(\mathcal{O};\mathbb{C}^n)},
\eqno(3.2)
$$
where $C_\mathcal{O}$ is the norm of the operator (3.1).
Let $S_\varepsilon$ be the smoothing operator (2.4).
By $R_\mathcal{O}$ we denote the operator of restriction of functions in $\mathbb{R}^d$
onto the domain $\mathcal{O}$. We put
$$
K_D(\varepsilon) = R_\mathcal{O} [\Lambda^\varepsilon] S_\varepsilon
b(\mathbf{D})P_\mathcal{O} (\mathcal{A}_{D}^0)^{-1}.
\eqno(3.3)
$$
The operator $b(\mathbf{D}) P_\mathcal{O} (\mathcal{A}_D^0)^{-1}$
is a continuous mapping of $L_2(\mathcal{O};\mathbb{C}^n)$ into $H^1(\mathbb{R}^d;\mathbb{C}^m)$.
Using Lemma 2.3 and relation $\Lambda \in \widetilde{H}^1(\Omega)$, it is easy to check that
the operator $[\Lambda^\varepsilon] S_\varepsilon$
is continuous from $H^1(\mathbb{R}^d;\mathbb{C}^m)$ to $H^1(\mathbb{R}^d;\mathbb{C}^n)$.
Hence, the operator (3.3) is continuous from $L_2(\mathcal{O};\mathbb{C}^n)$ to
$H^1(\mathcal{O};\mathbb{C}^n)$.

The following statement was proved in [PSu, (7.10)].

\smallskip\noindent\textbf{Proposition 3.1.} \textit{Let} $\mathcal{O} \subset \mathbb{R}^d$
\textit{be a bounded domain of class} $C^2$.
\textit{Let} $\mathbf{u}_\varepsilon$ \textit{be the solution of the problem} (1.4),
\textit{and let} $\mathbf{u}_0$ \textit{be the solution of the problem} (1.7)
\textit{with} $\mathbf{F} \in L_2(\mathcal{O};\mathbb{C}^n)$.
\textit{Let} $\widetilde{\mathbf{u}}_0 = P_\mathcal{O} \mathbf{u}_0$, \textit{where}
$P_\mathcal{O}$ \textit{is the extension operator} (3.1).
\textit{Let} ${\mathbf{w}}_\varepsilon \in H^1(\mathcal{O};\mathbb{C}^n)$
\textit{be the generalized solution of the problem}
$$
\mathcal{A}_{\varepsilon} \mathbf{w}_{\varepsilon} = 0 \ \text{in}\  \mathcal{O},
\quad \mathbf{w}_{\varepsilon} |_{\partial \mathcal{O}} =
\varepsilon \Lambda^{\varepsilon} S_\varepsilon b (\mathbf{D})
\widetilde{\mathbf{u}}_0 |_{\partial \mathcal{O}}.
\eqno(3.4)
$$
\textit{Then for} $0< \varepsilon \leq 1$ \textit{we have}
$$
\| \mathbf{u}_\varepsilon - \mathbf{u}_0 - \varepsilon
\Lambda^\varepsilon S_\varepsilon b(\mathbf{D})\widetilde{\mathbf{u}}_0 +
{\mathbf{w}}_\varepsilon\|_{H^1(\mathcal{O};\mathbb{C}^n)}
\leq \widetilde{C} \varepsilon \|\mathbf{F}\|_{L_2(\mathcal{O};\mathbb{C}^n)}.
\eqno(3.5)
$$
\textit{The constant} $\widetilde{C}$
\textit{depends only on} $m$, $d$, $\alpha_0$, $\alpha_1$, $\| g \|_{L_\infty}$,
$\| g^{-1} \|_{L_\infty}$, \textit{the parameters of the lattice} $\Gamma$
\textit{and the domain} $\mathcal{O}$.

\smallskip
The following theorem was proved in [PSu, Theorem 7.1].

\smallskip\noindent\textbf{Theorem 3.2.} \textit{Suppose that the assumptions of Theorem} 1.1
\textit{are satisfied. Let} $\widetilde{\mathbf{u}}_0 = P_\mathcal{O} \mathbf{u}_0$,
\textit{where} $P_\mathcal{O}$ \textit{is the extension operator} (3.1).
\textit{Then for} \hbox{$0< \varepsilon \leq \varepsilon_2$} \textit{we have}
$$
\| \mathbf{u}_{\varepsilon} - \mathbf{u}_0 - \varepsilon \Lambda^\varepsilon S_\varepsilon b(\mathbf{D}) \widetilde{\mathbf{u}}_0
\|_{H^1 (\mathcal{O}; \mathbb{C}^{n})} \leq C \varepsilon^{1/2}
\| \mathbf{F} \|_{L_2 (\mathcal{O}; \mathbb{C}^{n})},
\eqno(3.6)
$$
\textit{or, in operator terms,}
$$
\| \mathcal{A}_{D,\varepsilon}^{-1} - (\mathcal{A}^0_D)^{-1} -
\varepsilon K_D (\varepsilon) \|_{L_2 (\mathcal{O}; \mathbb{C}^{n}) \to
H^1 (\mathcal{O}; \mathbb{C}^{n})} \leq C {\varepsilon}^{1/2}.
$$
\textit{The constant} $C$
\textit{depends only on} $m$, $d$, $\alpha_0$, $\alpha_1$, $\| g \|_{L_\infty}$,
$\| g^{-1} \|_{L_\infty}$, \textit{the parameters of the lattice} $\Gamma$,
\textit{and the domain} $\mathcal{O}$.

\smallskip
Recall that $\left(\partial \mathcal{O} \right)_{\varepsilon}$
denotes the $\varepsilon$-neighborhood of $\partial \mathcal{O}$.
For sufficiently small $\varepsilon$, we fix two cut-off functions
$\theta_{\varepsilon} (\mathbf{x})$ and $\widetilde{\theta}_{\varepsilon} (\mathbf{x})$
in $\mathbb{R}^d$ such that
$$
\begin{aligned}
&\theta_{\varepsilon} \in C_0^{\infty} (\mathbb{R}^d), \quad
\text{supp}\,\theta_\varepsilon \subset (\partial \mathcal{O})_\varepsilon,
\quad 0 \leq \theta_{\varepsilon} (\mathbf{x}) \leq 1,
\\
&\theta_{\varepsilon} (\mathbf{x}) |_{\partial \mathcal{O}} = 1, \quad
\varepsilon \left| \nabla \theta_{\varepsilon} (\mathbf{x}) \right| \leq \kappa = \textrm{const};
\end{aligned}
\eqno(3.7)
$$
$$
\begin{aligned}
&\widetilde{\theta}_{\varepsilon} \in C_0^{\infty} (\mathbb{R}^d), \quad
\text{supp}\,\widetilde{\theta}_\varepsilon \subset (\partial \mathcal{O})_{2\varepsilon},
\quad 0 \leq \widetilde{\theta}_{\varepsilon} (\mathbf{x}) \leq 1,
\\
&\widetilde{\theta}_{\varepsilon} (\mathbf{x}) = 1\ \text{for}\
\mathbf{x} \in  (\partial \mathcal{O})_\varepsilon, \quad
\varepsilon \left| \nabla \widetilde{\theta}_{\varepsilon} (\mathbf{x}) \right| \leq \widetilde{\kappa} = \textrm{const}.
\end{aligned}
\eqno(3.8)
$$
We denote
$$
{\boldsymbol{\phi}}_\varepsilon = \varepsilon \theta_\varepsilon \Lambda^\varepsilon S_\varepsilon b(\mathbf{D}) \widetilde{\mathbf{u}}_0.
\eqno(3.9)
$$
From (1.2), (2.2), (2.6), (3.2), and (3.7) it follows that
$$
\| {\boldsymbol{\phi}}_\varepsilon\|_{L_2(\mathcal{O};\mathbb{C}^n)} \leq \varepsilon M \alpha_1^{1/2}
\| \widetilde{\mathbf{u}}_0 \|_{H^1(\mathbb{R}^d;\mathbb{C}^n)}
\leq \varepsilon M \alpha_1^{1/2} C_\mathcal{O} \widehat{c} \| \mathbf{F} \|_{L_2(\mathcal{O};\mathbb{C}^n)},
\eqno(3.10)
$$
cf. [PSu, (7.14)]. The norm of the function
(3.9) in $H^1(\mathcal{O};\mathbb{C}^n)$ was estimated in [PSu, Lemma 7.4].
A similar estimate is true if $\theta_\varepsilon$ is replaced by $\widetilde{\theta}_\varepsilon$.
We formulate the corresponding result.

\smallskip\noindent\textbf{Lemma 3.3.} \textit{Suppose that the assumptions of Theorem} 1.1
\textit{are satisfied. Let} $\theta_\varepsilon$ \textit{and}
$\widetilde{\theta}_\varepsilon$ \textit{be functions satisfying}
(3.7), (3.8), \textit{and let} ${\boldsymbol{\phi}}_\varepsilon$
\textit{be defined by} (3.9). \textit{Then we have}
$$
\| {\boldsymbol{\phi}}_\varepsilon\|_{H^1(\mathcal{O};\mathbb{C}^n)} \leq C_2
\varepsilon^{1/2} \| \mathbf{F} \|_{L_2(\mathcal{O};\mathbb{C}^n)},\quad 0< \varepsilon \leq \varepsilon_2,
\eqno(3.11)
$$
$$
\| \varepsilon \widetilde{\theta}_\varepsilon \Lambda^\varepsilon S_\varepsilon b(\mathbf{D}) \widetilde{\mathbf{u}}_0\|_{H^1(\mathcal{O};\mathbb{C}^n)} \leq
\widetilde{C}_2 \varepsilon^{1/2} \| \mathbf{F} \|_{L_2(\mathcal{O};\mathbb{C}^n)},\quad 0< 2\varepsilon \leq \varepsilon_2.
\eqno(3.12)
$$
\textit{The constants} $C_2$ \textit{and} $\widetilde{C}_2$
\textit{depend only on} $m$, $d$, $\alpha_0$, $\alpha_1$, $\| g \|_{L_\infty}$,
$\| g^{-1} \|_{L_\infty}$, \textit{the parameters of the lattice} $\Gamma$,
\textit{and the domain} $\mathcal{O}$.

\smallskip\noindent\textbf{3.2. Proof of Theorem 1.1. Step 1.}
Roughening (3.5), we obtain
$$
\| \mathbf{u}_\varepsilon - \mathbf{u}_0 - \varepsilon \Lambda^\varepsilon S_\varepsilon b(\mathbf{D}) \widetilde{\mathbf{u}}_0 + {\mathbf{w}}_\varepsilon\|_{L_2(\mathcal{O};\mathbb{C}^n)}
\leq \widetilde{C} \varepsilon \|\mathbf{F}\|_{L_2(\mathcal{O};\mathbb{C}^n)},\quad 0< \varepsilon \leq 1.
\eqno(3.13)
$$
Combining (1.2), (2.2), (2.6), and (3.2), we see that
$$
\|\Lambda^\varepsilon S_\varepsilon b(\mathbf{D})\widetilde{\mathbf{u}}_0\|_{L_2(\mathcal{O})} \leq M
\alpha_1^{1/2} C_\mathcal{O} \widehat{c} \|\mathbf{F}\|_{L_2(\mathcal{O})}.
\eqno(3.14)
$$
From (3.13) and (3.14) it follows that
$$
\| \mathbf{u}_\varepsilon - \mathbf{u}_0\|_{L_2(\mathcal{O})}
\leq \varepsilon(\widetilde{C} + M \alpha_1^{1/2} C_\mathcal{O} \widehat{c})
\|\mathbf{F}\|_{L_2(\mathcal{O})} + \| {\mathbf{w}}_\varepsilon\|_{L_2(\mathcal{O})},\quad 0< \varepsilon \leq 1.
\eqno(3.15)
$$
Therefore, the proof of estimate (1.8) is reduced to estimating
of ${\mathbf{w}}_\varepsilon$ in $L_2(\mathcal{O};\mathbb{C}^n)$.

For this purpose, we need the following lemma.

\smallskip\noindent\textbf{Lemma 3.4.}
\textit{Suppose that the assumptions of Theorem} 1.1 \textit{are satisfied. Let}
$\widetilde{\theta}_\varepsilon$
\textit{be a function satisfying} (3.8). \textit{Consider the operator}
$$
\widetilde{K}_D(\varepsilon) = R_\mathcal{O} [(1- \widetilde{\theta}_\varepsilon)\Lambda^\varepsilon]
S_\varepsilon b(\mathbf{D}) P_\mathcal{O} (\mathcal{A}^0_D)^{-1},
\eqno(3.16)
$$
\textit{which is a continuous mapping of} $L_2(\mathcal{O};\mathbb{C}^n)$ \textit{into}
$H^1_0(\mathcal{O};\mathbb{C}^n)$.
\textit{Let} $(\widetilde{K}_D(\varepsilon))^*: H^{-1}(\mathcal{O};\mathbb{C}^n)
\to L_2(\mathcal{O};\mathbb{C}^n)$
\textit{be the operator adjoint to the operator} (3.16), \textit{i.~e.},
$$
\begin{aligned}
\left((\widetilde{K}_D(\varepsilon))^* \mathbf{f}, \mathbf{v}\right)_{L_2(\mathcal{O})}=
\left( \mathbf{f}, \widetilde{K}_D(\varepsilon)\mathbf{v}\right)_{L_2(\mathcal{O})},
\\
\forall\,\mathbf{f}\in H^{-1}(\mathcal{O};\mathbb{C}^n),
\ \ \forall\,\mathbf{v} \in L_2(\mathcal{O};\mathbb{C}^n).
\end{aligned}
\eqno(3.17)
$$
\textit{Then the operator} $\mathcal{A}^{-1}_{D,\varepsilon}$,
\textit{viewed as a continuous mapping of} $H^{-1}(\mathcal{O};\mathbb{C}^n)$ \textit{into}
$L_2(\mathcal{O};\mathbb{C}^n)$, \textit{admits the following approximation}
$$
\| \mathcal{A}_{D,\varepsilon}^{-1} - (\mathcal{A}^0_D)^{-1} -
\varepsilon (\widetilde{K}_D(\varepsilon))^*
\|_{H^{-1} (\mathcal{O}) \to L_2 (\mathcal{O})} \leq (C+ \widetilde{C}_2) {\varepsilon}^{1/2},
\ \
0< 2\varepsilon \leq \varepsilon_2.
\eqno(3.18)
$$

\smallskip\noindent\textbf{Proof.}
From (3.6) and (3.12) it follows that
$$
\begin{aligned}
\| \mathbf{u}_{\varepsilon} - \mathbf{u}_0 -
\varepsilon (1 - \widetilde{\theta}_\varepsilon)\Lambda^\varepsilon
S_\varepsilon b(\mathbf{D})\widetilde{\mathbf{u}}_0
\|_{H^1 (\mathcal{O}; \mathbb{C}^{n})} \leq (C +\widetilde{C}_2) \varepsilon^{1/2}
\| \mathbf{F} \|_{L_2 (\mathcal{O}; \mathbb{C}^{n})},
\\
0< 2\varepsilon \leq \varepsilon_2.
\end{aligned}
\eqno(3.19)
$$
The function under the norm-sign on the left belongs to
$H^1_0(\mathcal{O};\mathbb{C}^n)$. In operator terms, (3.19) means that
$$
\| \mathcal{A}_{D,\varepsilon}^{-1} - (\mathcal{A}^0_D)^{-1} -
\varepsilon \widetilde{K}_D(\varepsilon)\|_{L_2 (\mathcal{O}) \to
H^1_0 (\mathcal{O})} \leq (C + \widetilde{C}_2) {\varepsilon}^{1/2},
\ \
0< 2\varepsilon \leq \varepsilon_2.
\eqno(3.20)
$$
This implies (3.18) by the duality arguments.
Indeed, combining (2.1), the similar identity for $(\mathcal{A}_D^0)^{-1}$
and (3.17), we see that for any $\mathbf{f} \in H^{-1}(\mathcal{O};\mathbb{C}^n)$ and
$\mathbf{v} \in L_2(\mathcal{O};\mathbb{C}^n)$ one has
$$
\begin{aligned}
&\left((\mathcal{A}_{D,\varepsilon}^{-1} - (\mathcal{A}^0_D)^{-1}
-\varepsilon (\widetilde{K}_D(\varepsilon))^*)\mathbf{f}, \mathbf{v}\right)_{L_2(\mathcal{O})}
\\
&=
\left(\mathbf{f}, (\mathcal{A}_{D,\varepsilon}^{-1} - (\mathcal{A}^0_D)^{-1}
-\varepsilon \widetilde{K}_D(\varepsilon)) \mathbf{v}\right)_{L_2(\mathcal{O})}.
\end{aligned}
$$
Together with (3.20) this yields
$$
\begin{aligned}
&\left| \left((\mathcal{A}_{D,\varepsilon}^{-1} - (\mathcal{A}^0_D)^{-1}
-\varepsilon (\widetilde{K}_D(\varepsilon))^*)\mathbf{f}, \mathbf{v}\right)_{L_2(\mathcal{O})} \right|
\\
&\leq (C + \widetilde{C}_2)\varepsilon^{1/2} \| \mathbf{f}\|_{H^{-1}(\mathcal{O})} \|\mathbf{v}\|_{L_2(\mathcal{O})},
\ \ \forall \, \mathbf{f}\in H^{-1}(\mathcal{O};\mathbb{C}^n),\ \ \forall \, \mathbf{v} \in L_2(\mathcal{O};\mathbb{C}^n),
\end{aligned}
$$
which implies (3.18). $\ \bullet$

\smallskip

From (3.7) and (3.9) it follows that
${\boldsymbol{\phi}}_\varepsilon \vert_{\partial \mathcal{O}} =
\varepsilon \Lambda^\varepsilon S_\varepsilon b(\mathbf{D})
\widetilde{\mathbf{u}}_0\vert_{\partial \mathcal{O}}$.
Then, by (3.4), the function
${\mathbf{w}}_\varepsilon - {\boldsymbol{\phi}}_\varepsilon$ is the solution of the problem
$$
\mathcal{A}_\varepsilon ({\mathbf{w}}_\varepsilon - {\boldsymbol{\phi}}_\varepsilon)
= {\mathbf{F}}_\varepsilon\  \text{in} \ \mathcal{O},\quad ({\mathbf{w}}_\varepsilon - {\boldsymbol{\phi}}_\varepsilon)\vert_{\partial \mathcal{O}} =0,
\eqno(3.21)
$$
where ${\mathbf{F}}_\varepsilon = - \mathcal{A}_\varepsilon {\boldsymbol{\phi}}_\varepsilon$.
It is easily seen that ${\mathbf{F}}_\varepsilon \in H^{-1}(\mathcal{O};\mathbb{C}^n)$, and
$$
\| {\mathbf{F}}_\varepsilon \|_{H^{-1}(\mathcal{O};\mathbb{C}^n)}
\leq \alpha_1 d^{1/2} \|g\|_{L_\infty}
\| {\boldsymbol{\phi}}_\varepsilon\|_{H^1(\mathcal{O};\mathbb{C}^n)},
\eqno(3.22)
$$
see [PSu, (4.15)]. From (3.11) and (3.22) it follows that
$$
\| {\mathbf{F}}_\varepsilon \|_{H^{-1}(\mathcal{O};\mathbb{C}^n)}
\leq C_3 \varepsilon^{1/2} \|\mathbf{F} \|_{L_2(\mathcal{O};\mathbb{C}^n)},
\quad 0< \varepsilon \leq \varepsilon_2,
\eqno(3.23)
$$
where $C_3 = \alpha_1 d^{1/2} \|g\|_{L_\infty} C_2$.
Note also that ${\mathbf{F}}_\varepsilon$ is supported in $(\partial \mathcal{O})_\varepsilon$.

Now we apply approximation (3.18) to the problem (3.21). Since
${\mathbf{w}}_\varepsilon - {\boldsymbol{\phi}}_\varepsilon =
\mathcal{A}_{D,\varepsilon}^{-1} {\mathbf{F}}_\varepsilon$, then
$$
\begin{aligned}
\| {\mathbf{w}}_\varepsilon - {\boldsymbol{\phi}}_\varepsilon -
(\mathcal{A}_{D}^0)^{-1} {\mathbf{F}}_\varepsilon
- \varepsilon (\widetilde{K}_D(\varepsilon))^* {\mathbf{F}}_\varepsilon\|_{L_2(\mathcal{O})}
\leq (C + \widetilde{C}_2) {\varepsilon}^{1/2} \| {\mathbf{F}}_\varepsilon\|_{H^{-1}(\mathcal{O})},
\\
0< 2\varepsilon \leq \varepsilon_2.
\end{aligned}
\eqno(3.24)
$$
By (3.16) and (3.17), for any $\mathbf{v}\in L_2(\mathcal{O};\mathbb{C}^n)$ we have
$$
\left((\widetilde{K}_D(\varepsilon))^* {\mathbf{F}}_\varepsilon,\mathbf{v}\right)_{L_2(\mathcal{O})}
=\left({\mathbf{F}}_\varepsilon, (1-\widetilde{\theta}_\varepsilon) \Lambda^\varepsilon S_\varepsilon
b(\mathbf{D}) P_\mathcal{O} (\mathcal{A}_D^0)^{-1}\mathbf{v}\right)_{L_2(\mathcal{O})}.
\eqno(3.25)
$$
Since $1 - \widetilde{\theta}_\varepsilon(\mathbf{x})=0$ for
$\text{dist}\,\{\mathbf{x},\partial \mathcal{O}\} \leq \varepsilon$,
and ${\mathbf{F}}_\varepsilon$ is supported in $(\partial \mathcal{O})_\varepsilon$,
then the right-hand side of (3.25)
is equal to zero. Consequently, $(\widetilde{K}_D(\varepsilon))^* {\mathbf{F}}_\varepsilon=0$.

Then (3.24) and (3.23) imply that
$$
\| {\mathbf{w}}_\varepsilon - {\boldsymbol{\phi}}_\varepsilon -
(\mathcal{A}_{D}^0)^{-1} {\mathbf{F}}_\varepsilon \|_{L_2(\mathcal{O})}
\leq (C + \widetilde{C}_2) C_3 {\varepsilon} \|\mathbf{F}\|_{L_2(\mathcal{O})},
\ \ 0< 2\varepsilon \leq \varepsilon_2.
\eqno(3.26)
$$
The norm of ${\boldsymbol{\phi}}_\varepsilon$ in $L_2(\mathcal{O};\mathbb{C}^n)$ admits estimate (3.10).
It remains to estimate the $L_2$-norm of the function $(\mathcal{A}_{D}^0)^{-1} {\mathbf{F}}_\varepsilon$.

\smallskip\noindent\textbf{3.3. Proof of Theorem 1.1. Step 2.}

\smallskip\noindent\textbf{Lemma 3.5.} \textit{Suppose that the assumptions of Theorem} 1.1
\textit{are satisfied. Let} ${\boldsymbol{\phi}}_\varepsilon$ \textit{be defined by} (3.9),
\textit{and let} ${\mathbf{F}}_\varepsilon = - \mathcal{A}_\varepsilon {\boldsymbol{\phi}}_\varepsilon =
- b(\mathbf{D})^* g^\varepsilon b(\mathbf{D}) {\boldsymbol{\phi}}_\varepsilon$.
\textit{Then the function} ${\boldsymbol{\eta}}_\varepsilon
:= (\mathcal{A}_D^0)^{-1} {\mathbf{F}}_\varepsilon$
\textit{satisfies the following estimate}:
$$
\| {\boldsymbol{\eta}}_\varepsilon\|_{L_2(\mathcal{O};\mathbb{C}^n)}
\leq C_4 \varepsilon \|\mathbf{F}\|_{L_2(\mathcal{O};\mathbb{C}^n)},
\quad 0< \varepsilon \leq \varepsilon_2.
\eqno(3.27)
$$
\textit{The constant} $C_4$ \textit{depends only on} $m$, $d$, $\alpha_0$, $\alpha_1$,
$\| g \|_{L_\infty}$, $\| g^{-1} \|_{L_\infty}$, \textit{the parameters of the lattice} $\Gamma$,
\textit{and the domain} $\mathcal{O}$.

\smallskip\noindent\textbf{Proof.} The function
${\boldsymbol{\eta}}_\varepsilon\in H^1_0(\mathcal{O};\mathbb{C}^n)$
is the generalized solution of the Dirichlet problem
$\mathcal{A}^0 {\boldsymbol{\eta}}_\varepsilon = {\mathbf{F}}_\varepsilon$,
${\boldsymbol{\eta}}_\varepsilon \vert_{\partial \mathcal{O}}=0$. It means that
$$
\begin{aligned}
\int_\mathcal{O} \langle g^0 b(\mathbf{D}) {\boldsymbol{\eta}}_\varepsilon, b(\mathbf{D}) \mathbf{h} \rangle
\,d\mathbf{x} = \int_\mathcal{O} \langle {\mathbf{F}}_\varepsilon, \mathbf{h} \rangle\,d\mathbf{x}=
-\int_\mathcal{O} \langle g^\varepsilon b(\mathbf{D}) {\boldsymbol{\phi}}_\varepsilon, b(\mathbf{D})
\mathbf{h} \rangle\,d\mathbf{x},
\\
\forall \,\mathbf{h} \in H^1_0(\mathcal{O};\mathbb{C}^n).
\end{aligned}
\eqno(3.28)
$$
If $\mathbf{h} \in H^2(\mathcal{O};\mathbb{C}^n) \cap H^1_0(\mathcal{O};\mathbb{C}^n)$,
then it is possible to integrate by parts in the left-hand side of (3.28). Hence,
$$
\begin{aligned}
\int_\mathcal{O} \langle {\boldsymbol{\eta}}_\varepsilon, \mathcal{A}^0 \mathbf{h} \rangle\,d\mathbf{x} =
-\int_\mathcal{O} \langle g^\varepsilon b(\mathbf{D}) {\boldsymbol{\phi}}_\varepsilon, b(\mathbf{D}) \mathbf{h} \rangle\,d\mathbf{x},
\\
\forall \,\mathbf{h} \in H^2(\mathcal{O};\mathbb{C}^n) \cap H^1_0(\mathcal{O};\mathbb{C}^n).
\end{aligned}
\eqno(3.29)
$$

Now we write down the norm of the function ${\boldsymbol{\eta}}_\varepsilon$ in
$L_2(\mathcal{O};\mathbb{C}^n)$ as the norm of continuous antilinear functional:
$$
\| {\boldsymbol{\eta}}_\varepsilon\|_{L_2(\mathcal{O})}
= \sup_{0 \ne \mathbf{G} \in L_2(\mathcal{O};\mathbb{C}^n)} \frac{\left| \int_\mathcal{O}
\langle {\boldsymbol{\eta}}_\varepsilon, \mathbf{G}\rangle \,d\mathbf{x} \right|}
{\|\mathbf{G}\|_{L_2(\mathcal{O})}}.
$$
We put $\mathbf{h}=(\mathcal{A}^0_D)^{-1} \mathbf{G}$, $\mathbf{G}\in L_2(\mathcal{O};\mathbb{C}^n)$.
Then $\mathbf{G} = \mathcal{A}^0 \mathbf{h}$, and $\mathbf{h}$ runs through
$H^2(\mathcal{O};\mathbb{C}^n) \cap H^1_0(\mathcal{O};\mathbb{C}^n)$ if
$\mathbf{G}$ runs through $L_2(\mathcal{O};\mathbb{C}^n)$ (see Subsection 2.2). Hence,
$$
\| {\boldsymbol{\eta}}_\varepsilon\|_{L_2(\mathcal{O})}
= \sup_{0 \ne \mathbf{h} \in H^2(\mathcal{O})\cap H^1_0(\mathcal{O})} \frac{\left| \int_\mathcal{O}
\langle {\boldsymbol{\eta}}_\varepsilon,
\mathcal{A}^0 \mathbf{h}\rangle \,d\mathbf{x} \right|}{\|\mathcal{A}^0\mathbf{h}\|_{L_2(\mathcal{O})}}.
\eqno(3.30)
$$
By (2.3), we have $\|\mathcal{A}^0\mathbf{h}\|_{L_2(\mathcal{O})} \geq (\widehat{c})^{-1}
\| \mathbf{h}\|_{H^2(\mathcal{O})}$. Combining this with (3.29) and (3.30), we obtain
$$
\| {\boldsymbol{\eta}}_\varepsilon \|_{L_2(\mathcal{O})} \leq
\widehat{c} \sup_{0 \ne \mathbf{h} \in H^2(\mathcal{O})\cap H^1_0(\mathcal{O})}
\frac{\left| \int_\mathcal{O} \langle g^\varepsilon b(\mathbf{D}) {\boldsymbol{\phi}}_\varepsilon,
b(\mathbf{D}) \mathbf{h}\rangle \,d\mathbf{x} \right|}{\| \mathbf{h} \|_{H^2(\mathcal{O})}}.
\eqno(3.31)
$$
Next, since $b(\mathbf{D})= \sum_{l=1}^d b_l D_l$ and, by (1.2), $|b_l| \leq \alpha_1^{1/2}$, then
$$
\|g^\varepsilon b(\mathbf{D}) {\boldsymbol{\phi}}_\varepsilon\|_{L_2(\mathcal{O})}
\leq \|g\|_{L_\infty} \alpha_1^{1/2} d^{1/2} \| {\boldsymbol{\phi}}_\varepsilon \|_{H^1(\mathcal{O})}.
\eqno(3.32)
$$
Taking into account that the function ${\boldsymbol{\phi}}_\varepsilon$ is supported in
the $\varepsilon$-neighborhood of $\partial \mathcal{O}$, from (3.31) and (3.32) we see that
$$
\| {\boldsymbol{\eta}}_\varepsilon \|_{L_2(\mathcal{O})} \leq
\widehat{c} \|g\|_{L_\infty} \alpha_1^{1/2} d^{1/2} \| {\boldsymbol{\phi}}_\varepsilon
\|_{H^1(\mathcal{O})} \sup_{0 \ne \mathbf{h} \in H^2(\mathcal{O})\cap H^1_0(\mathcal{O})}
\frac{\|b(\mathbf{D}) \mathbf{h}\|_{L_2(B_\varepsilon)}}{\| \mathbf{h} \|_{H^2(\mathcal{O})}}.
\eqno(3.33)
$$

Applying Lemma 2.2 and taking into account that
$|b(\mathbf{D}) \mathbf{h}| \leq \alpha_1^{1/2} \sum_{l=1}^d |D_l \mathbf{h}|$,
for $0< \varepsilon \leq \varepsilon_1$ we have:
$$
\begin{aligned}
&\int_{B_\varepsilon} |b(\mathbf{D}) \mathbf{h}|^2 \,d\mathbf{x}
\leq \alpha_1 d \sum_{l=1}^d \int_{B_\varepsilon} |D_l \mathbf{h}|^2\,d\mathbf{x}
\\
&\leq \alpha_1 d \beta \varepsilon \sum_{l=1}^d
\|D_l \mathbf{h}\|_{H^1(\mathcal{O})}\|D_l \mathbf{h}\|_{L_2(\mathcal{O})}
\leq \alpha_1 d \beta \varepsilon \| \mathbf{h}\|_{H^2(\mathcal{O})}
\| \mathbf{h}\|_{H^1(\mathcal{O})}.
\end{aligned}
$$
Hence,
$$
\sup_{0 \ne \mathbf{h} \in H^2(\mathcal{O})\cap H^1_0(\mathcal{O})}
\frac{\|b(\mathbf{D}) \mathbf{h}\|_{L_2(B_\varepsilon)}}{\| \mathbf{h} \|_{H^2(\mathcal{O})}}
\leq (\alpha_1 d \beta)^{1/2} \varepsilon^{1/2},\ \ 0< \varepsilon \leq \varepsilon_1.
\eqno(3.34)
$$
Finally, from (3.11), (3.33), and (3.34) it follows that (3.27) is valid with
$C_4 = \widehat{c} \|g\|_{L_\infty} \beta^{1/2} \alpha_1 d C_2$. $\ \bullet$

Now it is easy to complete the \textbf{proof of Theorem 1.1}.
By (3.10), (3.26), and (3.27), we have
$$
\| {\mathbf{w}}_\varepsilon \|_{L_2(\mathcal{O})} \leq C_5 \varepsilon \| \mathbf{F} \|_{L_2(\mathcal{O})},
\ \ 0 < 2\varepsilon \leq \varepsilon_2,
$$
where $C_5 = (C + \widetilde{C}_2) C_3 + M \alpha_1^{1/2} C_\mathcal{O} \widehat{c} + C_4$.
Combining this with (3.15), we arrive at (1.8) with
$C_1 = \widetilde{C} + M \alpha_1^{1/2} C_\mathcal{O} \widehat{c} + C_5$.
$\ \bullet$

\end{document}